\documentstyle{article}
\newtheorem{dfn}{Definition}[section]
\newtheorem{prop}[dfn]{Proposition}
\newtheorem{thm}[dfn]{Theorem}
\newtheorem{lem}[dfn]{Lemma}

\newtheorem{rem}[dfn]{Remark}

\newcommand{\bol}[1]{\mbox{\boldmath ${#1}$}}
\newcommand{\qed}{\hbox{\rule[-2pt]{3pt}{6pt}}}

\begin{document}
\title{Integral soluitons of $q$-difference equations \\
of the hypergeometric type with  $|q|=1$.
\thanks{To appear in the proceedings of the workshop 
 ``Infinite Analysis'' (Oct.15--19, 1996) at the IIAS, Japan } }
\author{Michitomo Nishizawa\thanks{694m5035@cfi.waseda.ac.jp}
$\mbox{ }$  and Kimio Ueno\thanks{uenoki@cfi.waseda.ac.jp}\\
Department of Mathematics,\\
School of Sciense and Engeneering,\\
Waseda University.}
\date{}
\maketitle

\begin{abstract}
 Two integral solutions of $q$-difference equations of the 
hypergeometric type with $|q|=1$ are constructed by using 
the double sine function. One is an integral of
the Barnes type and the other is of the Euler type.
\end{abstract}


\section{Introduction}
The hypergeometric $q$-difference equation is one of the most important
examples among $q$-difference systems and many studies have been achieved
\cite{gas}. However, these are concerned with the case that 
$0<q<1$. In the case of $|q|=1$, studies on $q$-difference
systems are not sufficiently explored. The difficulty comes from
 the facts that fundamental functions such as `` $q$-gamma function''
are not known in the case of $|q|=1 $.\par
 Recently, Jimbo and Miwa \cite{jim} have constructed an integral solution of 
the quantized Kniznik-Zamolodotikov equation with $|q|=1$. Inspired by the
result of Lukyanov \cite{luk}, they have given an integral solution
 by means of Kurokawa's double sine function \cite{kur}.
From a point of view of $q$-analysis, their work is very significant 
because it is thought of a first step of the study of $q$-difference
 system with $|q|=1$.\par
 In this article, we give two integral solutions of
$q$-difference equations of the hypergeometric type with $|q|=1$. 
One is an integral of the Barns type and the other is of the Euler type. 
 Once we obtain the $q$-gamma function with $|q|=1$,
 we can construct these integral representations in the same way
as in the case that $0<q<1$. Furthermore we can show that they are 
solutions of $q$-difference equations of the hypergeometric type
with $|q|=1$.\par 
 This article is organized as follows: In section 2, we give a survey of
integral representations of the hypergeometric series and the basic 
hypergeometric series with $0<q<1$. In section 3, we define the
 ``$q$-gamma function'' with $|q|=1$ by using the double sine function.
 In section 4, an integral of the Barnes type 
is introduced in the case of  $|q|=1$ and this function is shown to satisfy
the hypergeometric $q$-difference equation.
 In section 5, we consider an analogue of Euler's integral representation.
On this consideration, we must regard $q$-shifted factorials as the ``$q$-gamma
function'' with $|q|=1$, so it is needed to transform a multiplicative
variable to an additive variable. This integral gives a solution of the
differnce equation which is obtained by writing
 the hypergeometric $q$-difference equation by using an additive variable.\par 
 We would like to mention that our studies is significant when
 one considers the representation of the quantum group
 $SL_{q}(2,\bol{R})$. It is known that $q$ must be $|q|=1$
in $SL_{q}(2,\bol{R})$ (Masuda et.al. \cite{mas}), therefore, the
harmonic analysis on this quantum group should be closely linked to the
hypergeometric $q$-difference equation with $|q|=1$.  


\section{Preliminaries}
 In this section, we give a brief survey of integral representations
of the hypergeometric series
  \begin{equation}
       F(a,b,c;z)  =
      \sum_{k=0}^{\infty} \frac{(a)_{k}(b)_{k}}{(c)_{k} k!}z^{k}
     \qquad ( \mbox{for } |z|<1 ),
  \end{equation}
where $(a)_{k}:= a(a-1)\cdots(a-k+1)$,
and of the basic hypergeometric series with $0<q<1$
  \begin{equation}
    \phi(q^a,q^b,q^c;q,z)
     =\sum_{k=0}^{\infty}
       \frac{(q^a;q)_{k}(q^b;q)_{k}}{(q^c;q)_{k}(q;q)_{k}}z^k
       \qquad(\mbox{for }|z|<1),
  \end{equation}nnn
where $(a;q)_{k}:=\prod_{l=0}^{k}(1-aq^{l})$.
 
\subsection{Barnes' contour integral representation}
 Barnes' contour integral representation is so defined that  sum of residue
of the integrand is equal to the hypergeometric series.\par 
Let us define $(-z)^s:= \exp(s \log (-z))$, where we choose such a branch of 
logarithm that this logarithm takes real value when $z$ is 
on negative real line.
To define this integral, the following lemma is important.

\begin{lem}
$(1)$ The function $\pi(-z)^{s}/ \sin\pi s$ 
has simple poles at $s=k$ $( k\in \bol{Z})$, and the residue there
is $z^{k}$. \par
\noindent
$(2).$ We have, for 
$|z|<1$, that 
  \begin{equation}
    \frac{\pi (-z)^{s}}{\sin \pi s}
      =O\left[\exp\left\{-|\Im s|\arg (-z)\right\}\right]
   \end{equation}
as $\Im s \to \infty$ preserving $|\Re s| <\infty.$
\label{lem:21}\end{lem}

Let us fix a real number $\delta$ such that
$0<\delta<\pi$ and suppose $z$ to be in a sector $S_{1}:=\{z\in\bol{C}|
-\pi+\delta < \arg (-z) < \pi - \delta, |z|<1 \}.$  Barnes' contour 
integral of the hypergeometric series is given as follows:
  \begin{equation}
      F(a,b,c;z)= \frac{\Gamma(c)}{\Gamma(a)\Gamma(b)}
    \left(\frac{-1}{2\pi i}\right)
    \int_{-i\infty}^{i\infty}
      \frac{\Gamma(a+s)\Gamma(b+s)}{\Gamma(c+s)\Gamma(s+1)}
      \frac{\pi(-z)^{s}}{\sin \pi s} ds
\label{eqn:21}\end{equation}
where the contour lies on the right of poles
  \begin{equation}
     s=-an+n \qquad s=-b+n \qquad \left(n \in \bol{Z}_{\leq 0}\right)
  \end{equation}     
and on the left of poles $ s=m$ $(m\in \bol{Z}_{\geq 0}).$\par
Thanks to the Stirling formula of the gamma function and 
Lemma \ref{lem:21} (2), we can see that the integral (\ref{eqn:21}) 
converges uniformly in $S_1$ . 
Furthermore, by using deformation of the integral contour and
residue calculus based on Lemma \ref{lem:21} (1), one can show 
that the integral (\ref{eqn:21}) is the hypergeometric series. For the
details, see \cite{whi}. \par
Next, we consider a $q$-analogue of (\ref{eqn:21}) in the case that
$0<q<1.$  Let us put $q=e^{-2\pi\tau}$ $(\tau > 0).$ 
The counterpart of (4), which is known as Watson's contour
integral, is given as follows:
  \begin{equation}
      \phi(q^{a},q^{b},q^{c};q,z)=\frac{\Gamma(c;q)}{\Gamma(a;q)\Gamma(b;q)}
       \left(\frac{-1}{2\pi i}\right)
          \int_{-i\infty}^{i\infty}
       \frac{\Gamma(a+s;q)\Gamma(b+s;q)}{\Gamma(c+s;q)\Gamma(s+1;q)}
       \frac{\pi (-z)^{s}}{\sin \pi s} ds
  \label{eqn:24} \end{equation}
where $\Gamma(z;q)$ is the $q$-gamma function defined by
  \begin{equation}
    \Gamma(z:q):= \frac{(q;q)_{\infty}}{(q^{z};q)_{\infty}}
       (1-q)^{1-z},
  \end{equation}
and the contour lies on the right of poles
   $$ s=-a+n_{1}+\frac{n_{2}}{\tau}, \qquad
    s=-b+n_{1}+\frac{n_{2}}{\tau}, \qquad
      (n_{1} \in \bol{Z}_{\leq 0}, \quad n_{2} \in \bol{Z})$$
and on the left of poles $s=m$ $(m\in \bol{Z}_{\geq 0}).$\par
From Lemma \ref{lem:21} (2) and the fact that 
  $$ \left|\frac{\Gamma(a+s;q)\Gamma(b+s;q)}{\Gamma(c+s;q)\Gamma(1+s;q)}
   \right| \leq \mbox{Const.}
   \prod_{k=1}^{\infty}\frac{(1+e^{-(c+k+\Re s)\tau})(1+e^{-(1+k+\Re s)\tau})}
   {(1-e^{-(a+k+\Re s)\tau})(1-e^{-(b+k+\Re s)\tau})},$$
it follows that the integral (\ref{eqn:24}) converges uniformly
in $S_1$. By using the same technique, one can show that the integral
 (\ref{eqn:24}) is equal to the basic hypergeometric series \cite{gas}.

\subsection{Euler's integral representation}
 Euler's integral representation for the hypergeometric series is
  \begin{equation}
      F(a,b,c;z)
    =\frac{\Gamma(c)}{\Gamma(a)\Gamma(c-b)}
      \int_{0}^{1} t^{a-1}(1-t)^{c-a-1}(1-zt)^{-b}dt.
  \label{eqn:eul}\end{equation}
From the binomial theorem and an integral representation of the beta 
function, it follows that the integral (\ref{eqn:eul}) gives 
the hypergeometric series. \par
A $q$-analogue of this representation is given, by
using the Jackson integral, as follows:
  \begin{equation}
      \phi(q^a,q^b,q^c;q,z)
     = \frac{\Gamma(c;q)}{\Gamma(b;q)\Gamma(c-b;q)}
       \int_{0}^{1} t^{b}
       \frac{(tzq^{a};q)_{\infty}(tq;q)_{\infty}}
         {(tz;q)_{\infty}(tq^{c-b};q)_{\infty}}
       \frac{d_{q}t}{t}.
\label{eqn:qeul}\end{equation}
In the same way as the classical case, by using the $q$-binomial theorem
and the Jackson integral representation of the $q$-beta function,
we can prove that (\ref{eqn:qeul}) is equal to the basic hypergeometric 
series. 

\subsection{The hypergeometric $q$-difference equations}
 When $|q|<1$, the basic hypergeometric series $\phi(q^{a},q^{b},q^{c};q,z)$ 
is convergent for $|z|<1$, and satisfies the hypergeometric $q$-difference
equation;
\begin{equation}
  (L_{q}\phi)(z)=0,
\label{eqn:qhg}\end{equation}
where 
  \begin{eqnarray}
    & &[z]:=\frac{1-q^{z}}{1-q}, \qquad (T_{q}f)(z):=f(qz),\nonumber\\
    & &D_{q}:=\frac{1-T_{q}}{(1-q)z}, \qquad 
    [\vartheta + a]:= \frac{1-q^{a}T_{q}}{1-q}\nonumber\\
    & &L_{q}:= z^{-1}[\vartheta][\vartheta + c-1]
      -[\vartheta +a][\vartheta + b]\\[4pt]
    & &\qquad = z(q^{c}-q^{a+b+1}z)D_{q}^{2}\nonumber\\[4pt]
    & &\qquad  \quad- \left\{[c] - \frac{(1-q^{a})(1-q^{b})
      -(1-q^{a+b+1})}{1-q}z\right\}D_{q}
      - [a][b].\nonumber
  \end{eqnarray}
We should note that the basic hypergeometric series with $|q|=1$ is not
 convergent 
(so it gives only a formal solution to the hypergeometric $q$-difference
equation).


\section{``$q$-gamma funciton'' with $|q|=1$ }
 Let us define a function $\widetilde{\Gamma}(z;q)$ 
which satisfies 
  \begin{equation}
    \widetilde{\Gamma}(z+1;q)=[z]\widetilde{\Gamma}(z;q)
  \label{eqn:31}\end{equation}  
in the case of $|q|=1.$\par
 For this end, we need {\it the double zeta function}
$\zeta_{2}(s,z|\bol{\omega})$, {\it the double gamma function} 
$\Gamma_{2}(z|\bol{\omega})$ and {\it the double sine function} 
$S_{2}(z|\bol{\omega})$  (cf. \cite{bar}, \cite{jim}, \cite{kur}, \cite{shi}).

\begin{dfn}
For $\bol{\omega}:=(\omega_{1}, \omega_{2}) \in \bol{C}^{2}$, we define 
$\zeta_{2}(s,z|\bol{\omega})$, $\Gamma_{2}(z|\bol{\omega})$ and 
$S_{2}(z|\bol{\omega})$ by
  \begin{eqnarray*}
   & & \zeta_{2}(s,z|\bol{\omega})
         :=\sum_{m_{1},m_{2}\in \bol{Z}_{\geq 0}}
         \left(z+m_{1}\omega_{1}+m_{2}\omega_{2}\right)^{-s}, \\
   & & \Gamma_{2}(z|\bol{\omega})
         := \exp \left( \frac{\partial}{\partial s}
            \zeta_{2}(s,z|\bol{\omega})\vert_{s=0} \right),\\
   & & S_{2}(z|\bol{\omega})
        := \Gamma_{2}(z|\bol{\omega})^{-1} 
           \Gamma_{2}(\omega_{1}+\omega_{2}-z|\bol{\omega}).
  \end{eqnarray*}
\end{dfn}
 
It is known that the double sine function satisfies the functional
relation
  \begin{equation}
    \frac{S_{2}(z+\omega_{1}|\bol{\omega})}{S_{2}(z|\bol{\omega})}
      =\frac{1}{2 \sin \frac{\pi z}{\omega_{2}}}.
  \end{equation}
Thus, we can construct a function satisfying (\ref{eqn:31}) by using
$S_{2}(z|\bol{\omega})$. We suppose that $|q|=1$ and that $q$ is not
a root of unity. Let us put $q=e^{2\pi i \omega}$ $(0<\omega<1, 
\omega \notin \bol{Q}).$ 

\begin{dfn}
 We set
  \begin{equation}
     \widetilde{\Gamma}(z;q):= (q-1)^{1-z}i^{z-1}q^{\frac{z(z-1)}{4}}
     S_{2}(z|(1, \frac{1}{\omega}))^{-1},
  \end{equation}
\end{dfn}

\noindent
which has the following properties.

\begin{prop}
$(1)$ $\widetilde{\Gamma}(z;q)$ has simple zeros
 at $z=n_{1}+\frac{n_{2}}{\omega} 
\quad (n_{1}, n_{2}\in \bol{Z}_{>0})$, and has simple poles at  
$z=n_{1}+\frac{n_{2}}{\omega} \quad (n_{1},n_{2}\in \bol{Z}_{\leq 0}).$\\
\noindent
$(2)$ $\widetilde{\Gamma}(z;q)$ satisfies the functional relation 
(\ref{eqn:31}).\\
\noindent
$(3)$ If we take $z\to \infty$ as $z$ is in any sector not
containing real line then \par $\widetilde{\Gamma}(z;q)$ has the following 
asymptotic behavior.
\begin{eqnarray*}
  & &\widetilde{\Gamma}(z;q) = \exp \left[ (1-z)\log(q-1) +(z-1)\log i
    \right. \\[4pt]
  & &\left.\qquad  \qquad +\frac{z(z-1)}{4}\log q 
    \mp \pi i  \left\{ \frac{\omega z^{2}}{2}-\frac{\omega+1}{2}z
    \right\} +O(1)\right] \quad (for \pm \Im z >0).
\end{eqnarray*}
\end{prop}
 
\noindent
These properties follow from the facts in the papers \cite{jim}, \cite{shi}.

\begin{rem}
In the case that \ $0<q<1$, we can also define $\widetilde{\Gamma}(z,q)$
by Definition 3.2 (in this case, $\omega=it, \quad t>0$).
Of course, we can see that $\widetilde{\Gamma}(z,q)
=C(z,q)\Gamma(z,q)$, where $C(z,q)$ is a function satisfying 
$C(z+1,q)=C(z,q)$ (cf. \cite{shi}). 
\end{rem}


\section{An integral representation of the Barnes type with $|q|=1$}

\subsection{Definition of the integral}
In order that the integral makes sense, we impose some conditions on
parameters $a,b$ and $c$.\par
\vspace{10pt}
\noindent
{\bf Conditions on the parameters } {\it $(B1)$ If we define sets $A_{1}$ 
and $A_{2}$ of the parameters by}
  $A_{1}:=\{a,b\}$, $A_{2}:=\{c,1\},$
{\it then we suppose that }
  $$''\Re \alpha > \Re \beta \qquad \mbox{for} \quad \forall\alpha \in A_{1},
    \quad \forall\beta \in A_{2}''$$
{\it or} 
  $$''\Im \alpha \ne \Im \beta \qquad \mbox{for} \quad \forall\alpha \in A_{1},
    \quad \forall\beta \in A_{2}''$$\par
\noindent
$(B2)$ {\it We suppose that} $$\omega \Re (a+b-c+1) < 1.$$ 
\par    
Under these conditions, we can define a function $\Phi(a,b,c;q,z)$
 in the same way as Bernes' contour integral by using
 $\widetilde{\Gamma}(z,q)$.

\begin{dfn}
Let us fix a real number $\delta$ such that $0<\delta<\pi-\pi\omega
\Re (a+b-c+1).$
For $z$ in the sector $S :=\{z\in\bol{C}| -\pi+\delta < \arg (-z) < \pi-
2\pi\omega \Re (a+b-c+1), |z|<1\},$  we define $\varphi(a,b,c;q;s,z)$ and
$\Phi(a,b,c;q,z)$ by
  \begin{eqnarray}
    & &\varphi(a,b,c;q;s,z) 
         := \frac{\widetilde{\Gamma}(a+s;q)\widetilde{\Gamma}(b+s;q)}
                 {\widetilde{\Gamma}(c+s;q)\widetilde{\Gamma}(1+s;q)}
            \frac{\pi(-z)^{s}}{\sin\pi s}\nonumber\\
    & &\Phi(a,b,c;q,z)
         := \frac{\widetilde{\Gamma}(c;q)}
                 {\widetilde{\Gamma}(a;q)\widetilde{\Gamma}(b;q)}
            \left(\frac{-1}{2\pi i}\right)
            \int_{-i\infty}^{i\infty}
            \varphi(a,b,c:q;s,z)ds
  \end{eqnarray}
where the contour lies on the right of poles
  $$ s=-a+n_{1}+\frac{n_{2}}{\omega}, \quad
     s=-b+n_{1}+\frac{n_{2}}{\omega}  \quad
     (n_{1}, n_{2} \in \bol{Z}_{\leq 0}) $$
and on the left of poles 
  \begin{eqnarray*}
    & & s=-c+n_{1}+\frac{n_{2}}{\omega}, \quad
        s=-1+n_{1}+\frac{n_{2}}{\omega}  \quad
        (n_{1}, n_{2} \in \bol{Z}_{> 0}), \\
    & & s=m \quad (m \in \bol{Z}_{\geq 0}).
  \end{eqnarray*}
\label{dfn:41}\end{dfn}

By using Lemma \ref{lem:21} (2) and Proposition 3.3 (2), it is shown that
  $$ \varphi(a,b,c;q;s,z) = O \left[ \exp(-\delta |s|)\right] 
     \quad \mbox{as} \quad s\to\pm i \infty$$
under the condition (B1). Thus the integral (15) converges uniformly
in $S$, and the analytic continuation (also
denote $\Phi(a,b,c;q,z)$) defines a many-valued analytic function of $z$.   

\subsection{The hypergeometric $q$-difference equation with $|q|=1$}
We prove that $\Phi(a,b,c;q,z)$ is a solution of the hypergeometric 
$q$-difference equation with $|q|=1$. We also use the notation (10) 
in the case that $|q|=1$.  

\begin{thm}
$\Phi(z):=\Phi(a,b,c;q,z)$ satisfies the hypergeometric $q$-difference
equation
  $$(L_{q}\Phi)(z)=0.$$
\end{thm}
{\it   Outline of Proof:  }
From the condition (B2), it follows that the action of $L_{q}$  commute 
with the integration. On the other hand, straightfoward calculation 
shows that the integrand $\varphi(a,b,c;q;s,z)$ satisfies 
the relation
\begin{equation}
  (L_{q}\varphi)(a,b,c;q;s,z)
      =\varphi(a+1,b+1,c;q;s-1,z)
        -\varphi(a+1,b+1,c;q;s,z).
\end{equation}
By means of Cauchy's theorem, one can verify that the integral of the
right-hand side of (16) vanishes. 
\qed


\section{An integral representation of the Euler type with $|q|=1$}
\subsection{Definition of the integral}
First let us recall Euler's integral (9) in the case that $0<q<1$. 
If we transform the variables $z$ and $t$ in the integrand of (\ref{eqn:qeul})
to $q^{x}$ and $q^{s}$ respectively, then we have
  $$(\mbox{the integrand of } (9))
    = \mbox{ Const. } 
      \frac{\Gamma(s+x;q)\Gamma(s+c-b;q)}
      {\Gamma(s+x+a;q)\Gamma(s+1;q)}q^{bs}.$$
Therefore, in the case that $|q|=1$, we consider the integral
  \begin{equation}
    \int \frac{\widetilde{\Gamma}(s+x;q)\widetilde{\Gamma}(s+c-b;q)}
       {\widetilde{\Gamma}(s+x+a;q)\widetilde{\Gamma}(s+1;q)}
       q^{bs}ds
  \label{eqn:51}\end{equation}
as a counterpart of (9). In order that the integral makes sense,
 we impose the following conditions on the parameters $a$, $b$ and $c$.\par
\vspace{10pt}
\noindent
{\bf Conditions on the parameters.} \ $(E1)$ $b-c \notin \bol{R}_{>0}, $ 
\ $(E2)$ $a \notin \bol{R}_{<0},$ \\ $(E3)$ $\Re b >0,$ $\Re(a-c-1)>0.$\par	
\vspace{10pt}
\noindent
Under these conditions, we can take such a suitable contour that 
the integral (\ref{eqn:51}) makes sense.

\begin{dfn}
For $x\notin\bol{R}_{<0}$, we define a function $\Psi(a,b,c;q,x)$ by
  \begin{equation}
    \Psi(a,b,c;q,x):=
       \int_{-i\infty}^{i\infty}
       \frac{\widetilde{\Gamma}(s+x;q)\widetilde{\Gamma}(s+c-b;q)}
       {\widetilde{\Gamma}(s+x+a;q)\widetilde{\Gamma}(s+1;q)}
       q^{bs}ds
  \end{equation}
where the contour lies on the right of the poles 
  $$ s=-x+n_{1}+\frac{n_{2}}{\omega}, \quad
     s=b-c+n_{1}+\frac{n_{2}}{\omega}, \quad
     (n_{1}, n_{2} \in \bol{Z}_{\leq 0}),$$
and on the left of the poles
  $$ s=-x-a+n_{1}+\frac{n_{2}}{\omega}, \quad
     s=-1+n_{1}+\frac{n_{2}}{\omega}, \quad
     (n_{1},n_{2} \in \bol{Z}_{>0}).$$
\end{dfn}

Thanks to the conditions (E1) and (E2), we can take the contour
of Definition 5.1. From the condition (E3), it follows that the
integral (18) converges uniformly and defines a single-valued
analytic function of $x$. 
       
\subsection{The difference equation for $\Psi(a,b,c;q,x)$}
Let us present an equation which $\Psi(a,b,c;q,x)$ satisfies.
For this end, we write the hypergeometric 
$q$-difference equation by using the ``additive'' variable $x$. We employ
the following notatiaons;
  \begin{eqnarray*}
    & & (T_{+}g)(x):=g(x+1), \quad
        [\vartheta+a]_{+}:= \frac{1-q^{a}T_{+}}{1-q}\\[5pt]
    & & L_{+}:= q^{-x}[\vartheta]_{+}[\vartheta+c-1]_{+}
                 -[\vartheta+a]_{+}[\vartheta+b]_{+}\\[4pt]
    & & \qquad   = \frac{1}{(1-q)^{2}}
                \left[ (q^{c-1-x}-q^{a+b})
                \left\{T_{+}^{2}-(1+q)T_{+}+q\right\}\right.\\
    & & \qquad \qquad  - \left\{(1-q^{c})q^{-x}+(1-q^{a})(1-q^{b})
                 -(1-q^{a+b+1})\right\}(T_{+}-1)\\
    & & \qquad \qquad  \left.- (1-q^{a})(1-q^{b})\right]  
  \end{eqnarray*}
Then the next theorem holds.

\begin{thm}
$\Psi(x):= \Psi(a,b,c;q,x)$ satisfies the difference equation
  $$(L_{+}\Psi)(x) = 0.$$
\end{thm}
This theorem can be proved just like Theorem 5.1.

\end{document}